\input amstex
\documentstyle{amsppt}
\magnification 1200
\vcorrection{-9mm}
\input epsf
\NoBlackBoxes

\topmatter
\title      Signatures of iterated torus links
\endtitle
\author     S.~Yu.~Orevkov
\endauthor
\abstract
We compute the multivariate signatures of any Seifert link (that is a union of some fibers
in a Seifert homology sphere), in particular, of the union of a torus link
with one or both of its cores (cored torus link).
The signatures of cored torus links are used in Degtyarev-Florens-Lecuona
splicing formula for computation of multivariate signatures of cables over links.
We use Neumann's computation of equivariant signatures of such links. 

For signatures of torus links with the core(s)
we also rewrite the Neumann's formula in terms of integral points
in a certain parallelogram, similar to Hirzebruch's
formula for signatures of torus links (without cores)
via integral points in a rectangle.
\endabstract

\address
Steklov Mathematical Institute, Gubkina 8, Moscow, Russia
\endaddress

\address
IMT, l'universit\'e Paul Sabatier, 118 route de Narbonne, Toulouse, France
\endaddress

\address
AGHA Laboratory, Moscow Institute of Physics and Technology, Russia
\endaddress

\endtopmatter

\def\Z{\Bbb Z}
\def\Q{\Bbb Q}
\def\R{\Bbb R}
\def\C{\Bbb C}
\def\sph{\Bbb S}
\def\ld{{(\hskip-.1em(}}
\def\rd{{)\hskip-.1em)}}

\def\card{\#} \def\Card{\card}
\def\li{\lfloor}
\def\ri{\rfloor}
\def\eps{\varepsilon}
\def\Null{\operatorname{n}}

\def\one{\bold 1}
\def\Re{\operatorname{Re}}
\def\Im{\operatorname{Im}}
\def\lk{\operatorname{lk}}
\def\trans{T}
\def\sign{\operatorname{sign}}
\def\Log{\operatorname{Log}}
\def\ind{\operatorname{ind}}

\def\seif{{\frak{Seif}}}  

\document

\def\refBri   {1}
\def\refCim   {2}
\def\refCF    {3}
\def\refCo    {4}
\def\refCNT   {5}
\def\refDFL   {6}
\def\refDFLii {7}
\def\refDFLiii{8}
\def\refEN    {9}
\def\refFlo   {10}
\def\refGo    {11}
\def\refLit   {12}
\def\refMat   {13}
\def\refNeuE  {14}
\def\refNeuI  {15}
\def\refOre   {16}
\def\refViro  {17}

\def\figHopf           {1}
\def\figSeif           {2}
\def\figDefect         {3}
\def\figDede           {4}
\def\figParallelogramm {5}
\def\figSplicePQ       {6}
\def\figSpliceAA       {7}
\def\figSpliceDE       {8}
\def\figSpliceFG       {9}
\def\figExSplice       {10}
\def\figExParal        {11}
\def\figExSign         {12}
\def\figExSignRes      {13}
\def\figExSignDiag     {14}
\def\figExTwo          {15}
\def\figExTwoDiag      {16}

\def\sectSplice  {2}
\def\sectIter    {3}
\def\sectDefSig  {4}
\def\sectDFL     {5}
\def\sectDefESig {6}
\def\sectDede    {7}
\def\sectNeumann {8}
\def\sectHirz    {9}
\def\sectSeif    {10}
\def\sectExample {11}

\def\defSplice    {\sectSplice.1}
\def\thEN         {\sectSplice.2}
\def\propLN       {\sectSplice.3}

\def\remColor     {\sectDefSig.1}
\def\propSigSym   {\sectDefSig.2} 
\def\propSigRev   {\sectDefSig.3} 
\def\propSigSig   {\sectDefSig.4} 
\def\propSigAlx   {\sectDefSig.5} 
\def\corSigAlx    {\sectDefSig.6} 
\def\lemSR        {\sectDefSig.7} \let\lemSimpleRoot=\lemSR
\def\corSR        {\sectDefSig.8} \let\corSimpleRoot=\corSR

\def\thDFL        {\sectDFL.1}

\def\propMatumoto {\sectDefESig.1}
\def\remOrient    {\sectDefESig.2} 

\def\lemDede      {\sectDede.1}

\def\thNeumann    {\sectNeumann.1}
\def\corNeumann   {\sectNeumann.2}
\def\remNeumann   {\sectNeumann.3}

\def\propHirz     {\sectHirz.1}
\def\remHirz      {\sectHirz.2}
\def\remDeepNest  {\sectHirz.3}
\def\remMultilink {\sectHirz.4}

\def\propColor   {\sectSeif.1}

\def\eqDefSeif    {1}
\def\eqSym        {2}
\def\eqDefAprime  {3}
\def\eqAlex       {4}
\def\eqPropHirz   {5}

\def\eqAlex       {7}
\def\eqCondA      {8}
\def\eqMM         {9}
\def\eqSpliceD    {10}
\def\eqSpliceE    {11}
\def\eqSpliceF    {12}
\def\eqSpliceG    {13}
\def\eqCondB      {14}
\def\eqSpliceFF   {15}
\def\eqSpliceGG   {16}

\def\eqExample    {17}
\def\eqDFL        {18}

\head 1. Introduction \endhead

An {\it iterated torus link} in the 3-sphere $\sph^3$ is a link obtained from the unknot\
by successive cabling operations either with the core removed or retained
(see~\S\sectIter\ for a precise definition).
The most important examples are links of singularities of plane complex analytic curves,
and links at infinity of plane affine algebraic curves.

If $K$ is an iterated torus knot, the 
Levine-Tristram signatures $\sigma_\zeta(K)$
(defined as the signature of $(1-\zeta)V+(1-\bar\zeta)V^\trans$ for a Seifert matrix
$V$) can be recursively computed using:
\roster
\item"$\bullet$"
Hirzebruch's formula (see [\refBri; \S6], [\refMat; \S4])
which expresses the signatures
of a torus knot $T(p,q)$ in terms of integral points in the rectangle $p\times q$;
\smallskip
\item"$\bullet$"
Litherland's formula [\refLit; Thm.~2] which expresses the signatures of a
$(p,q)$-cable of any knot $K$ via those of $K$ and $T(p,q)$
\endroster

To extend this approach to iterated torus links, it is natural to consider a multivariate
generalization of the Levine-Tristram signatures discussed
in [\refCF, \refDFL, \refFlo, \refViro].
In [\refDFL, \refDFLii],
an analog of the Litherland's formula is obtained for cables over links.
It expresses the multivariate signatures of a cable over $L$ via those
of $L$ and those of the union of a torus link with one or both of its cores
(the link $T_{m_1,m_2}(p,q)$, $m_1,m_2\in\{0,1\}$,
in the notation of \S\sectHirz\ below), but the latter signatures have been unknown.

The initial aim of this paper was to fill this gap and to compute
the multivariate signatures of $T_{m_1,m_2}(p,q)$.
This (and a little more) is done in Proposition \propColor\ which
reduces the multivariate signatures of any Seifert link
(in particular of $T_{m_1,m_2}(p,q)$)
to Tristram-Levine signatures of this link. The latter were already computed by
Walter Neumann in [\refNeuE], [\refNeuI].

More precisely,
In Proposition \propColor, for any positive Seifert
$n$-component link $L$ (that is a union of $n$ positively oriented fibers in a Seifert homology sphere), we show that its multivariate signature $\sigma_L$ 
considered as a function on the open cube ${]}0,1{[}^n$, is constant
on each member of a family of parallel hyperplanes transverse to the main diagonal,
so, all values of $\sigma_L$ are determined by the values on the main diagonal,
which coincide with the Levine-Tristram signatures up to a certain additive constant.
The proof of Proposition \propColor\ is based on the splice formula for multivariate signatures (which is the main result of [\refDFL]) combined with the observation that
a Seifert link can be spliced in many different ways.

In Proposition \propHirz\
we show that in the case of a torus link $T_{m_1,m_2}(p,q)$, the Neumann's formula
can be rewritten in such a way that it becomes almost identical with the Hirzebruch's formula but the parallelogram spanned by
$(p+m_1,m_2)$ and $(m_1,q+m_2)$ is used instead of the $p\times q$ rectangle.

In \S\S2--8 we give necessary definitions and results from
[\refDFL--\refDFLiii, \refEN, \refMat--\refNeuI].
In \S\sectHirz\ we prove Proposition~\propHirz;
in \S\sectSeif\ we prove Proposition~\propColor;
in \S\sectExample\ we consider some examples
and formulate some questions for further research.

\smallskip
\subhead Acknowledgement
\endsubhead
I am grateful to Alexander Degtyarev and Vincent Florens for useful discussions.


\head\sectSplice. Splicing and splice diagrams (after Eisenbud and Neumann)
\endhead

Let $(\Sigma',K'\cup L')$ and $(\Sigma'',K''\cup L'')$ be two links in homology spheres.
Here $K'$ and $K''$ are components of the respective links.
Let $T(K')$
and $T(K'')$ be their tubular neighbourhood disjoint from $L'$
and $L''$ respectively.

\medskip\noindent{\bf Definition \defSplice.}
The {\it splice} of $(\Sigma',K'\cup L')$ and $(\Sigma'',K''\cup L'')$
along $K'$ and $K''$ is introduced in [\refEN] (see also [\refNeuI], [\refDFL])
as $(\Sigma,L'\cup L'')$ where
$\Sigma=(\Sigma'\setminus T(K'))\cup_\varphi(\Sigma''\setminus T(K''))$
and $\varphi:T(K')\to T(K'')$ is a homeomorphism which identifies
the meridian of $K'$ with the longitude of $K''$ and vice versa.
\medskip

\midinsert
  \centerline{\epsfxsize=35mm\epsfbox{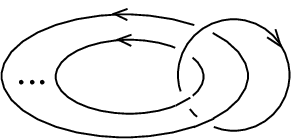}}
\botcaption{Figure \figHopf} $\seif(0,1,\dots,1)$ (called in [\refDFL]
      a generalized Hopf link)
\endcaption
\endinsert

Let $a_1,\dots,a_n$ be positive pairwise coprime integers.
Following [\refEN, \refNeuI],
we define {\it Seifert link} $\seif(a_1,\dots,a_n)$ to be the link
$(\Sigma,S_1\cup\dots\cup S_n)$ where $\Sigma=\Sigma(a_1,\dots,a_n)$
is the unique Seifert fibered homological $3$-sphere which has fibers
$S_1,\dots,S_n$ of degrees $a_1,\dots,a_n$ and no other fibers of degree $>1$
(see [\refEN, \S7] for a detailed construction).
The orientation of $\Sigma$ is chosen so that
all the linking numbers $\lk(S_i,S_j)$, $i\ne j$, are positive.
This definition extends to the case when
$a_1,\dots,a_n$ are any pairwise coprime integers by setting
$$
\split
  &\Sigma(a_1,\dots,-a_i,\dots,a_n)=-\Sigma,\\
  &\seif(a_1,\dots,-a_i,\dots,a_n)=
(-\Sigma,S_1\cup\dots\cup(-S_i)\cup\dots\cup S_n),
\endsplit
$$
and $\seif(0,1,\dots,1)$ to be the link in the $3$-sphere shown in
Figure~\figHopf. Finally, we define
$$
   \seif(a_1,\dots,a_k; a_{k+1},\dots,a_n)=(\Sigma,S_1\cup\dots\cup S_k),
                                                                    \eqno(\eqDefSeif)
$$
that is $\seif(a_1,\dots,a_n)$ with the last $n-k$ components removed.

\midinsert
\centerline{\epsfxsize=25mm\epsfbox{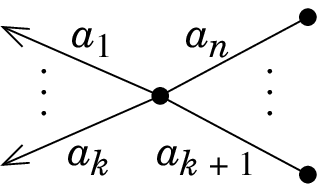}}
\botcaption{Figure \figSeif}
  The splice diagram of $\seif(a_1,\dots,a_k;a_{k+1},\dots,a_n)$.
\endcaption
\endinsert

As shown in [\refEN], any graph link in a {\it graph homology sphere}
(i.e., a graph manifold which is a homology sphere)
can be obtained as the result of splicing together Seifert links.
The way of splicing is encoded in [\refEN]
by certain decorated trees (or forests) called {\it splice diagrams}.
The splice diagram of $\eps\seif(a_1,\dots,a_k;a_{k+1},\dots,a_n)$, $\eps=\pm1$, is
the graph in Figure~\figSeif\ if $\eps=+1$.
If $\eps=-1$, then we use the white color for the central vertex.
If the orientation of some component is reversed, we put the minus sign
near the corresponding arrowhead vertex (if the sign is not shown,
we assume that it is plus).
The splicing of $(\Sigma',L'\cup K')$
with $(\Sigma'',L''\cup K'')$ along $K'$ and $K''$ has the splice
diagram $\boxed{\Gamma'}\!{-}\!{-}\!{-}\!\boxed{\Gamma''}$ where
$\boxed{\Gamma'}\!{\longrightarrow}$ and
${\longleftarrow}\!\boxed{\Gamma''}$ are the splice diagrams of the links
being spliced and the indicated arrowheads correspond to $K'$ and $K''$
respectively.

Let $-\Gamma$ denote the diagram $\Gamma$ with the opposite signs
of all the arrowhead vertices. 

\proclaim{ Theorem \thEN } {\rm(See [\refEN,~Thm.~8.1].)}
Two splice diagrams determine the same link if and only
if they are obtained from each other by a sequence of the following equivalences:
\roster
\item
      $\Gamma\approx -\Gamma$,
\vskip10pt
\item
     \hskip-30pt
     \lower8pt\hbox to 96pt{\vbox{\noindent
         ${\boxed{\Gamma}}\!{>}\!\!{\bullet}\!{-}\!{-}\!{-}\!{-}\!%
          {\boxed{\Gamma'}}$\par\vskip-6pt\noindent\hskip26pt$a$%
     }}
     $\approx$
     \hskip-30pt
     \lower8pt\hbox to 120pt{\vbox{\noindent
         ${\boxed{\Gamma}}\!{>}\!{\circ}\!{-}\!{-}\!{-}\!{-}\!%
          {\boxed{-\Gamma'}}$\par\vskip-6pt\noindent\hskip23pt$-a$%
     }}
\vskip10pt
\item
      $\;\boxed{\Gamma}\!{>}\!\!{\bullet}\!{-}\!{-}\!{-}\!{-}\!{-}\!{-}\!{\bullet}$
      $\;\approx\;$
      $\boxed{\Gamma}\!{>}\!\!{\bullet}$
      \par\vskip-5pt\noindent\hskip30pt$1$
\vskip10pt
%
     \hskip-43pt
     \lower8pt\hbox to 120pt{\vbox{%
         \vskip-8pt\noindent
         ${\boxed{\Gamma'}}\!{-}\!{-}\!{-}\!{-}\!{\bullet}\!{-}\!{-}\!{-}\!{-}\!%
          {\boxed{\Gamma''}}$\par\vskip-4pt\noindent\hskip24pt$a'$%
           \hskip18pt$a''$%
     }}
     $\approx\;$
         ${\boxed{\Gamma'}}\!{-}\!{-}\!{-}\!{-}\!{-}\!{-}\!%
          {\boxed{\Gamma''}}$
\vskip10pt
\item
     \hskip-30pt
     \lower20pt\hbox to 90pt{\vbox{%
        \noindent$\boxed{\Gamma_1}$\par\vskip-1pt\noindent
           \hskip12pt $1\setminus$\par\vskip-11pt
        \noindent\hskip3pt${\vdots}$\hskip13pt${\bullet}\!{-}\!{-}\!{-}\!{-}\!%
             {\bullet}$\par\vskip-5pt\noindent
            \hskip12pt $1/$\hskip10pt$0$\par\vskip-1pt
        \noindent$\boxed{\Gamma_n}$
     }}
     $\approx$
     \hskip-30pt
     \lower20pt\hbox to 70pt{\vbox{%
        \noindent${\boxed{\Gamma_1}}\!{-}\!{-}\!{-}\!{\bullet}$\par
        \noindent\hskip13pt${\vdots}$\par\vskip3pt
        \noindent${\boxed{\Gamma_n}}\!{-}\!{-}\!{-}\!{\bullet}$\par
     }}
     (disjoint union)
\vskip10pt
\item
     $\;\boxed{\Gamma}\;$ $\approx$ $\;\boxed{\Gamma}\;\;$
     ${\bullet}\!{-}\!{-}\!{-}\!{-}\!{\bullet}$ (disjoint union);
\vskip10pt
\item
     \hskip-30pt
     \lower20pt\hbox to 140pt{\vbox{%
        \noindent$\boxed{\Gamma_1}$\hskip59pt\boxed{\Gamma_{k+1}}
        \par\vskip-2pt\noindent
           \hskip8pt $a_1\!{\setminus}\hskip51pt{/}a_{k+1}$
        \par\vskip-12pt
        \noindent\hskip3pt${\vdots}$\hskip13pt${\bullet}\!{-}\!{-}\!{-}\!{-}\!%
                \!{-}\!{-}\!{-}\!{-}\!{\bullet}$\hskip27pt${\vdots}$
             \par\vskip-5pt\noindent
            \hskip6pt $a_k/$\hskip8pt$b_1$\hskip18pt$b_2$\hskip7pt$\setminus a_n$%
       \par\vskip-0.8pt
        \noindent$\boxed{\Gamma_k\hskip-1pt}$\hskip59pt$\boxed{\,\;\Gamma_n\;\,}$
     }}
     $\approx$
     \hskip-30pt
     \lower20pt\hbox to 92pt{\vbox{%
        \noindent$\boxed{\Gamma_1}$\hskip7.5pt\boxed{\Gamma_{k+1}}
        \par\vskip-2pt\noindent
           \hskip8pt $a_1\!{\setminus}{/}a_{k+1}$
        \par\vskip-12pt
        \noindent\hskip3pt${\vdots}$\hskip13pt${\bullet}$\hskip27pt${\vdots}$
             \par\vskip-5pt\noindent
            \hskip6pt $a_k{/}{\setminus}a_n$%
       \par\vskip-0.8pt
        \noindent$\boxed{\Gamma_k\hskip-1pt}$\hskip8.5pt$\boxed{\,\;\Gamma_n\;\,}$
     }}
     if $\displaystyle\prod_{i\le k} a_i=b_2$ and
        $\displaystyle\prod_{i>k} a_i=b_1$.
\endroster
\endproclaim
\medskip

\proclaim{ Proposition \propLN } {\rm(See [\refEN, \S10].)}
Let $(\Sigma,S_1\cup\dots\cup S_k)$ be as in
(\eqDefSeif). Then the linking numbers are
$\lk(S_i,S_j)=a_1\dots\hat a_i\dots\hat a_j\dots a_n$, $1\le i<j\le k$
(as usual, the hat means the omission of the corresponding factor).
\endproclaim

In fact, in [\refEN, \S10], the linking numbers are expressed in terms of
a splice diagram for any graph link in a graph homology $3$-sphere.


\head\sectIter. Cabling and iterated torus links  \endhead

Given an oriented link $L\cup K$ in a homology sphere $\Sigma$
where $K$ is a knot and a pair of integers $(p,q)\ne(0,0)$, we set
$d=\gcd(p,q)$, $p_1=p/d$, $q_1=q/d$, and
we define the {\it $(p,q)$-cabling of $L\cup K$ along $K$ with the core removed} (resp.
{\it with the core retained}) as $L\cup L_{p,q}$ (resp. $L\cup K\cup L_{p,q}$) where
$L_{p.q}$ is the union of $d$ disjoint knots $K_1\cup\dots\cup K_d$
such that for some tubular neighbourhood $T$ of $K$ disjoint from $L$, one has
$K_j\subset\partial T$, $[K_j]=p_1[K]\in H_1(T)$, and $\lk(K_j,K)=q_1$
for each $j=1,\dots,d$.

An iterated torus link in $\sph^3$ is defined as a link obtained from the unknot
by successive cabling. Note that reversion of orientation of a component is
equivalent to the $(-1,0)$-cabling with the core removed.

The $(p,q)$-cabling of $(\Sigma,L\cup K)$ along $K$ with the core retained
(resp. removed) is equivalent to the splicing of $(\Sigma,L\cup K)$ with
$$
    \seif(q_1,\underset{d}\to{\underbrace{1,\dots,1}},p_1)\quad
\text{(resp.~with}\;\;
    \seif(q_1,\underset{d}\to{\underbrace{1,\dots,1}};p_1))
$$
along $(K,S_1)$ where $S_1$ is the component corresponding to the weight $q_1$
(see [\refEN, Prop.~9.1] and the paragraph after it).


\head\sectDefSig. Multivariate link signatures
\endhead

Let $L$ be an oriented link in a homology sphere $\Sigma$ and let $V$ be the Seifert
form on a {\sl connected} Seifert surface. Let $\zeta\in\C\setminus\{1\}$, $|\zeta|=1$.
Then the Levine-Tristram signature $\sigma_\zeta(L)$ and nullity $\Null_\zeta(L)$ are
defined as the signature and nullity of $(1-\zeta)V + (1-\bar\zeta)V^\trans$.

Let $\sph^1$ be the unit circle in $\C$ and let
$\Cal T=\{e^{2\pi i\theta}\mid \theta\in\Q$ and
$0<\theta<1\}\subset\sph^1\setminus\{1\}$.
Let $L_1,\dots,L_\mu$ be the components of $L$.
A {\it multivariate signature and nullity} of $L$ defined in [\refCF, \refDFL, \refFlo] are functions $\sigma_L,\Null_L : \Cal T^\mu\to\Z$. As shown in [\refViro]
 (see also [\refCNT]),
these functions can be defined on $(\sph^1\setminus\{1\})^\mu$.
It is natural to extend them to $(\sph^1)^\mu$ by interpreting the value $1$
of the $i$-th argument as the removal of $L_i$ (see [\refDFL, \refDFLii]).

\medskip\noindent{\bf Remark \remColor. }
In fact, the multivatiate signatures are defined in [\refCF, \refFlo, \refViro]
for any {\it colored link}, i.e., a link $L$ with a fixed decomposition into
a disjoint union of sublinks (not necessarily connected) $L=L_1\cup\dots\cup L_\mu$.
In this case the $i$-th argument of $\sigma_L$ and $\Null_L$ corresponds to $L_i$.
In this paper we consider only the case when each $L_i$ is connected.
However, everything can be easily extended to the case of arbitrary colored links
due to [\refCF, Prop.~2.5].
\medskip

We refer to [\refViro], [\refCNT], or [\refDFLiii] for a definition of
$\sigma_L$ and $\Null_L$. Here we just mention some properties of them.

\proclaim{ Proposition \propSigSym } {\rm(See [\refCF, Prop.~2.7].)}
$$
\matrix
   \sigma_L(u_1,\dots,u_\mu) = \sigma_L(u_1^{-1},\dots,u_\mu^{-1}),\\
    \Null_L(u_1,\dots,u_\mu) =  \Null_L(u_1^{-1},\dots,u_\mu^{-1}).
\endmatrix                                                           \eqno(\eqSym)
$$
\endproclaim

\proclaim{ Proposition \propSigRev } {\rm(See [\refCF, Prop.~2.8].)}
Let $L'$ be obtained from
$L=L_!\cup\dots\cup L_\mu$ by reversing the orientation of $L_i$. Then
$\sigma_{L}(u)=\sigma_{L'}(u')$ and
$ \Null_{L}(u)= \Null_{L'}(u')$ where $u=(u_1,\dots,u_\mu)$
and $u'=(u_1,\dots,u_i^{-1},\dots,u_\mu)$.
\endproclaim

\proclaim{ Proposition \propSigSig } {\rm(See [\refCF, Prop.~2.5].)}
Let $\lambda\in\sph^1\setminus\{1\}$. Then
$$
   \sigma_\lambda(L) = \sigma_L(\lambda,\dots,\lambda)
                     - \sum_{1\le i<j\le\mu}\lk(L_i,L_j),
   \qquad
   \Null_\lambda(L) = \Null_L(\lambda,\dots,\lambda).
$$
\endproclaim

The following fact is proven in [\refCF] for links in the $3$-sphere only but
the proof extends (with certain efforts) to links in any homology sphere.

\proclaim{ Proposition \propSigAlx } {\rm(See [\refCF].)}
There exists a matrix $A=A(t_1,\dots,t_\mu)$ whose entries are Laurent
polynomials and such that:
\roster
\item"$\bullet$" the multivariate Alexander polynomial $\Delta_L(t_1,\dots,t_n)$
   is equal to $\det A$ up to some factors of the form $\pm t_i^{\pm1}$ or $\pm(t_i-1)$;
\item"$\bullet$" for any $u\in (\sph^1\setminus\{1\})^\mu$, the signature and nullity
   of $A(u)$ are $\sigma_L(u)$ and $\Null_L(u)$ respectively.
\endroster
\endproclaim

\proclaim{ Corollary \corSigAlx }
If the multivariate Alexander polynomial $\Delta_L(t_1,\dots,t_\mu)$ is
not identically zero, then $\sigma_L$ is constant on each connected component
of the complement of the zero set of $\Delta_L$ in
$(\sph^1\setminus\{1\})^\mu$.
\endproclaim

The following fact should be known and it is easy to prove.

\proclaim{ Lemma \lemSimpleRoot }
Let $A(t)=(a_{ij}(t))$, $t\in\R$, be a Hermitian matrix such that
$\Re a_{ij}$ and $\Im a_{ij}$ are real analytic functions of $t$.
Suppose that $\det A(t)$ has a simple root at $t=t_0$, and there are
no other roots in an interval $(t_0-2\eps,t_0+2\eps)$.
Let $s_0$ and $s_\pm$ be the signature of $A(t)$ at $t=t_0$ and $t=t_0\pm\eps$
respectively. Then $|s_+ - s_-|=2$, $s_0=(s_++s_-)/2$, and the nullity of $A(t_0)$ is $1$.
\qed
\endproclaim

By combining Proposition~\propSigAlx\ with Lemma~\lemSR, we obtain:

\proclaim{ Corollary \corSR } Let $u=(u_1,\dots,u_\mu)\in(\sph^1\setminus\{1\})^\mu$
be such that the gradient of $\Delta_L(t_1,\dots,t_\mu)$ at $u$ is non-zero.
Then, for a small neighbourhood $U\subset(\sph^1\setminus\{1\})^\mu$ of $u$, the
restrictions $\sigma_L|_U$ and $\Null_L|_U$ depend only on the sign of $\Delta_L$,
and one has $|s_1 - s_{-1}|=2$, $s_0=(s_1 + s_{-1})/2$, $n_{\pm1}=0$, and $n_0=1$
where $s_t=\sigma_L|_{U_t}$ and $n_t=\Null_L|_{U_t}$ for
$U_t=U\cap\{\sign\Delta_L=t\}$, $t=-1,0,1$.
\endproclaim


\head\sectDFL. Multivariate signatures of a splice
               (after Degtyarev, Florens and Lecuona)
\endhead

For $\ell=(\ell_1,\dots,\ell_\mu)\in\Z^\mu$, we define the {\it defect function}
$\delta_\ell:(\sph^1)^\mu\to\Z$ by setting
$$
   \delta_\ell(u_1,\dots,u_\mu) = \ind\Big(\sum\ell_i\Log u_i\Big)
       - \sum\ell_i\ind(\Log u_i)
$$
(see Figure~\figDefect) where $\Log:\sph^1\to[0,1{[}$ and $\ind:\R\to\Z$ are defined by
$$
   \Log(e^{2\pi i t})=t, \qquad \ind(x)=\li x\ri - \li-x\ri.
$$

\midinsert
\centerline{\epsfxsize=90mm\epsfbox{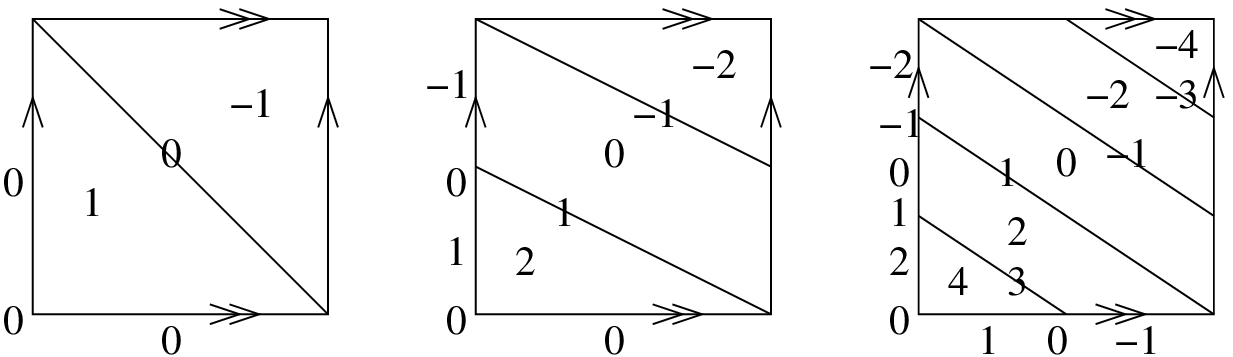}}
\centerline{$\delta_{1,1}(u)$ \hskip20mm
            $\delta_{1,2}(u)$ \hskip20mm
            $\delta_{2,3}(u)$}
\botcaption{Figure \figDefect} The values of $\delta_\ell$ on $(\sph^1)^2$
   for $\ell=(1,1)$, $(1,2)$, and $(2,3)$.
\endcaption
\endinsert

Let $(\Sigma,L'\cup L'')$ be the splice of $(\Sigma',K'\cup L')$ and
$(\Sigma'',K''\cup L'')$ along $K'$ and $K''$ (see Definition~\defSplice).
Let $L'_1,\dots,L'_{\mu'}$ and $L''_1,\dots,L''_{\mu''}$
be the components of $L'$ and $L''$ respectively.
Let $\ell'$ and $\ell''$ be the vectors of linking numbers
$$
  \ell'=\big(\lk(K',L'_1),\dots,\lk(K',L'_{\mu'})\big),\qquad
  \ell''=\big(\lk(K'',L''_1),\dots,\lk(K'',L''_{\mu''})\big).
$$

\proclaim{ Theorem \thDFL } (a). {\rm(See [\refDFL, Thm.~2.2], [\refDFLii, Thm.~5.2].)}
Assume that $L'$ and $L''$ are non-empty.
Let $u'\in(\sph^1)^{\mu'}$ and $u''\in(\sph^1)^{\mu''}$ be such that
$(v',v'')\ne(1,1)$ where $v'=(u')^{\ell'}$ and $v''=(u'')^{\ell''}$. Then
$$
\split
   &\sigma_{L'\cup L''}(u',u'') = \sigma_{K'\cup L'}(v'',u')+
      \sigma_{K''\cup L''}(v',u'')+\delta_{\ell'}(u')\delta_{\ell''}(u''),\\
   &\Null_{L'\cup L''}(u',u'') = \Null_{K'\cup L'}(v'',u')+
      \Null_{K''\cup L''}(v',u'').
\endsplit
$$

(b). {\rm(See [\refDFL, Addendum~2.7].)} Assume that $L'=\varnothing$ and
$L''\ne\varnothing$. Then, for any $u\in(\sph^1)^{\mu''}$ one has
$$
\split
    &\sigma_{(\Sigma,L'')}(u) = \sigma_{K'}(u^{\lambda''})
               + \sigma_{K''\cup L''}(u),\\
    &\Null_{(\Sigma,L'')}(u) = \Null_{K'}(u^{\lambda''})
               + \Null_{K''\cup L''}(u).
\endsplit
$$
\endproclaim

In the case when $(v',v'')=(1,1)$, Theorem~\thDFL(a) does not apply but
the signature and nullity of the splice
can be however often computed by the same formulas with a correction term which
ranges in $[-2,2]$. The correction term is a function of the values of
a certain invariant called {\it slope} evaluated on each of the splice components,
see details in [\refDFLii, \refDFLiii] (see also Question~\sectExample.1 below).


\head\sectDefESig. Equivariant link signatures
\endhead

Another kind of link signatures considered in this paper are
{\it equivariant signatures} $\sigma_\lambda^\pm$. We define them for fibered links
only though the definition can be extended to the general case
(see surveys [\refCo], [\refGo]).
So, let $L$ be a fibered link, $F$ be the fiber (in particular, $\partial F=L$),
and $h:H_1(F)\to H_1(F)$ be the monodromy operator. Then the one-variable Alexander
polynomial $\Delta_L(t)$ is the characteristic polynomial of $h$.
Let $H=H_1(F)\otimes\C$ and let $H=\bigoplus_\lambda H_\lambda$
be the splitting of $H$ according to the eigenvalues of $h$.
Let $l$ be the Seifert form on $H_1(F)$ extended to a sesquilinear form on $H$.
Then $\sigma_\lambda^+$ (resp. $\sigma_\lambda^-$) is defined as the signature of
the hermitian form $l+l^*$ (resp. $i(l-l^*)$) restricted to $H_\lambda$.

\proclaim{ Proposition \propMatumoto } {\rm(See [\refMat], [\refNeuE].)}
For any fibered link $L$, we have:

\smallskip
(a). $\sigma_\lambda^+ =  \sigma_{\bar\lambda}^+$ and
$\sigma_\lambda^-=\sigma_\lambda^+\sign\Im\lambda$,
in particular
     $\sigma_\lambda^- = -\sigma_{\bar\lambda}^-$ and
 $\sigma_{\pm1}^-=0$;

\smallskip
(b). if $|\lambda|\ne 1$ or $\lambda$ is not a root of the Alexander polynomial,
then $\sigma_\lambda^+=\sigma_\lambda^-=0$;

\smallskip
(c). Let $\omega = e^{i\varphi}$, $0<\varphi\le\pi$.
If $h$ is semisimple or $\omega$ is not
a root of $\Delta_L(t)$, then the Levine--Tristram signature and
nullity are
$$
    \sigma_\omega(L) = \sigma_\omega^+
    + 2\sum_{\smallmatrix 0\le\theta<\varphi\\\lambda=e^{i\theta}
      \endsmallmatrix}\sigma_\lambda^+
    \qquad\text{ and }\qquad \Null_\omega(L) = \dim H_\omega,
$$
i.e., $\Null_\omega(L)$ is the multiplicity of $\omega$ as a root
of $\Delta_L(t)$.
\endproclaim

\medskip\noindent{\bf Remark \remOrient. }
The orientation conventions in [\refMat] and those in [\refNeuE, \refNeuI]
are different. In [\refMat] they are chosen so that $l^*h=l$
whereas in [\refNeuE, \refNeuI] so that $lh=l^*$. By this reason,
the signs of $\sigma_\lambda^-$ in these sources also differ.
We use the convention from [\refMat] which seems to be more common nowadays.


\head\sectDede. The sawtooth function and an identity for it
\endhead
Let
$$
    \ld x \rd = \cases {1\over2} - x + \li x\ri, &x\not\in\Z\\
                        0,                &x\in\Z \endcases
$$
where $\li x\ri=\max\{n\in\Z\mid n\le x\}$ is the integer part of $x$.
Let, for a set $\Omega$,
$$
       \one_\Omega(x) = \cases 1, & x\in\Omega,\\ 0, & x\not\in\Omega. \endcases
$$
\proclaim{ Lemma \lemDede } Let $x,y,z_1,\dots,z_n\in\R$. Suppose that
$x+y+z_1+\dots +z_n=s\in\Z$ and $0<z_i<1$ for all $i=1,\dots,n$. Then
$$
    \ld x\rd+\ld y\rd+\ld z_1\rd+\dots+\ld z_n\rd = n/2 - \card\big(\Z\cap{]}x,s-y{[}\,\big)
        - \tfrac12\one_{\Z}(x) - \tfrac12\one_{\Z}(y).
$$
\endproclaim

\midinsert
\vbox{
  \centerline{$z_1+\dots+z_n$\hskip 10mm}\vskip-1mm
  \centerline{\epsfxsize=110mm\epsfbox{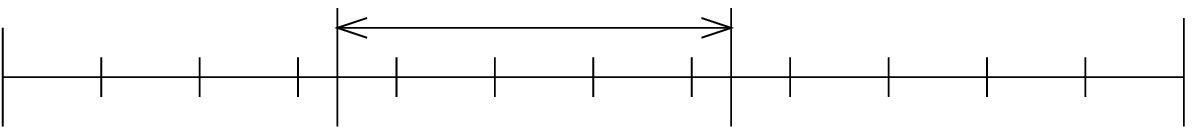}}
  \vskip-3mm
  \centerline{$0$\hskip 7mm
    \lower-3mm\hbox{$1$}\hskip 13.5mm
    \lower-3mm\hbox{$\li x\ri$}\hskip 1mm $x$\hskip 33mm $s{-}y$\hskip 0mm
    \lower-3mm\hbox{$s{-}\li y\ri$}\hskip 28mm $s$}}
\botcaption{Figure \figDede}
\endcaption
\endinsert

\demo{ Proof }
We consider the case $x,y\not\in\Z$ and we leave the other cases to the reader.
Since $\li z_1\ri=\dots=\li z_n\ri=0$, the left hand side is equal to
$$ \frac{n+2}{2} - (x - \li x\ri) - (y - \li y\ri) - \sum z_i
    = \frac{n+2}{2} - s + \li x\ri + \li y\ri
$$
and we have
$1+\card\big(\Z\cap{]}x,s-y{[}\,\big) = 
   (s-\li y\ri) - \li x\ri$ (see Figure \figDede).
\qed\enddemo


\head\sectNeumann. Neumann's formula for equivariant signatures
 of Seifert links
\endhead

Let $L=\seif(a_1,\dots,a_k;a_{k+1},\dots,a_n)$ with positive $a_1,\dots,a_n$.
We set:
$$
    m_j = \cases 1, &j\le k,\\ 0, &j>k, \endcases\qquad
    a'_j = (a_1\dots a_n)/a_j, \qquad
    m = \sum_{j=1}^n m_j a'_j.                        \eqno(\eqDefAprime)
$$

Let us choose $b_1,\dots,b_n$ so that
$b_j a'_j \equiv 1 \mod a_j$ for each $j=1,\dots,n$, and let
$$
   s_j = (m_j - b_j m)/a_j, \qquad j=1,\dots,n.
$$
\proclaim{ Theorem \thNeumann } {\rm(Neumann [\refNeuE, \refNeuI].)}
The link $L$ is fibered, its monodromy operator is semisimple, and one has:
$\sigma_1^+(L) = 1-k$, $\sigma_{-1}^+(L)=0$, and, for
$\lambda\ne 1$,
$$
   \sigma_\lambda^-(L) = \cases
       -2\sum_{j=1}^n \ld s_jk/m\rd
                    &\text{ if $\lambda = \exp(2\pi i k/m)$ with $k\in\Z$,}\\
       0            &\text{ otherwise}
   \endcases
$$
where $\ld\dots\rd$ is the sawtooth function defined in \S\sectDede.
\endproclaim

Recall that the sign of $\sigma_\lambda^-$ in [\refNeuE, \refNeuI] is opposite,
see Remark~\remOrient.

\proclaim{ Corollary \corNeumann }
For $\omega=e^{i\varphi}$, $0<\varphi\le\pi$, one has
$$
    \sigma_\omega(L) = 1-k + \sigma_\omega^-(L) +
    2\sum_{\smallmatrix0<\theta<\varphi\\\lambda=e^{i\theta}\endsmallmatrix}
    \sigma_\lambda^-(L)
$$
and $\Null_\omega(L)$ is equal to the multiplicity
of $\omega$ as a root of $\Delta_L(t)$.
\endproclaim

By [\refEN, Thm.~12.1], if $k\ge2$, the multivariate Alexander polynomial
of $L$ is
$$
  \Delta_L(t_1,\dots,t_k) = \big(t_1^{a'_1}\dots t_k^{a'_k}-1\big)^{n-2}
  \prod_{j=k+1}^n\big(t_1^{a'_1/a_j}\dots t_k^{a'_k/a_j}-1\big)^{-1}     \eqno(\eqAlex)
$$
and hence the one-variable Alexander polynomial is
$\Delta_L(t)=(t-1)\Delta_L(t,\dots,t)$. (If $k=1$, then the right hand side of (\eqAlex)
should be multiplied by $t_1-1$.)
The formula (\eqAlex) is a specialization of a general formula given in [\refEN]
for the Alexander polynomial of any graph multilink in a graph homology sphere.

\medskip\noindent{\bf Remark \remNeumann. }
In fact, the results of [\refNeuI, Thms.~5.1--5.3] are much stronger than
Theorem~\thNeumann: for any algebraic link (not only for positive Seifert links)
a decomposition of the Hermitian isometric structure into irreducible factors is
computed there, which includes a description of the Seifert form up to congruence
over $\C$ (there is a misprint in
[\refNeuI, Thm.~5.3]: the factor $i$ should be omitted in
the Seifert part of $\Lambda_\lambda^2$).


\head\sectHirz.
     Hirzebruch-type formula for signatures of a torus link with the core(s)
\endhead

Given four integers $a,b,c,d$ with $ad-bc > 0$, we denote the {\bf open} parallelogram in $\R^2$
spanned by the vectors $(a,b)$ and $(c,d)$ by $\Pi=\Pi(a,b,c,d)$. Thus
$$
    \Pi(a,b,c,d) = \{s(a,b)+t(c,d)\mid 0<s<1 \text{ and } 0<t<1\}.
$$
Let $u:\R^2\to\R$ be the linear function such that $u(a,b)=u(c,d)=1$.
For $0\le\theta<1$, we set (see Figure \figParallelogramm):
$$
\split
    N_\theta^- &= N_\theta^-(a,b,c,d) =
       \card\{(x,y)\in\Z^2\cap\Pi\mid \theta < u(x,y) < \theta+1\},\\
    N_\theta^+ &= N_\theta^+(a,b,c,d) =
       \card\{(x,y)\in\Z^2\cap\Pi\mid u(x,y) < \theta \text{ or }\theta+1 < u(x,y)\},\\
    M_\theta^- &= M_\theta^-(a,b,c,d) =
       \card\{(x,y)\in\Z^2\cap\Pi\mid u(x,y) = \theta+1\},\\
    M_\theta^+ &= M_\theta^+(a,b,c,d) =
       \card\{(x,y)\in\Z^2\cap\Pi\mid u(x,y) = \theta\}.
\endsplit
$$

\midinsert
\centerline{\epsfxsize=75mm\epsfbox{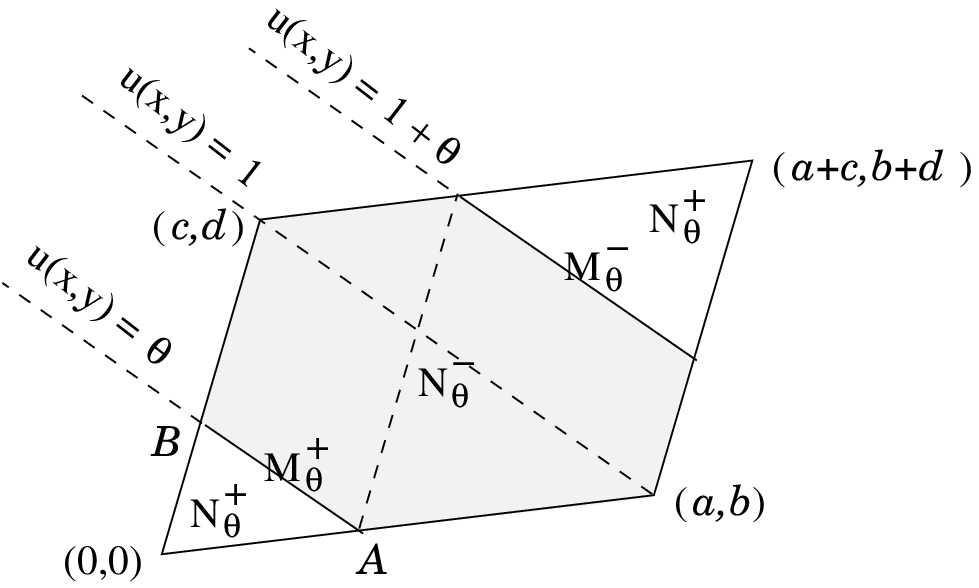}}
\botcaption{Figure \figParallelogramm} Definition of $N^{\pm}_\theta$ and $M_\theta^\pm$
\endcaption
\endinsert

Let $p$ and $q$ be positive integers (not necessarily coprime) and let
$m_1,m_2\in\{0,1\}$.
We define the torus link $T_{m_1,m_2}(p,q)$ as the intersection of the unit sphere in
$\C^2$ with the algebraic curve $\{(z_1,z_2)\mid z_1^{m_1}z_2^{m_2}(z_1^p - z_2^q)=0\}$
endowed with the boundary orientation induced from the intersection of the curve
with the unit ball, thus the linking number of any two components is positive.
If $m_1=m_2=0$, this is the torus link $T(p,q)$.

\proclaim{ Proposition \propHirz }  Let
$(a,b) = (p,0)+(m_1,m_2)$ and $(c,d) = (0,q) + (m_1,m_2)$. 
Then the equivariant signature and Levine-Tristram signature and nullity of
$T_{m_1,m_2}(p,q)$ at $\lambda=e^{2\pi i\theta}$ for $0<\theta<1$ are:
$$
   \sigma_\lambda^-= M^+_\theta - M^-_\theta,             \qquad
   \sigma_\lambda  = N_\theta^+ - N_\theta^- - m_1 - m_2, \qquad
    \Null_\lambda  = M_\theta^+ + M_\theta^-
$$
where $N_\theta^\pm = N_\theta^\pm(a,b,c,d)$, $M_\theta^\pm = M_\theta^\pm(a,b,c,d)$.
\endproclaim

\noindent{\bf Remark \remHirz.}
In the case $m_1=m_2=0$ this is Hirzebruch's formula
(see [\refBri; \S6] for $\theta=1/2$ and [\refMat; \S4] for any $\theta$).

\smallskip
\noindent{\bf Remark \remDeepNest.}
In [\refOre; Prop.~8.1], I gave the following formulas for the signature
and nullity of the braid closure $L$ of the braid $\Delta^n$ with $2k+1$ strings:
$$
    \sigma(L) = \sigma_{-1}(L) =
        \cases -nk(k+1) + (-1)^{(n-1)/2} &\text{if $k\equiv n\equiv 1 \mod 2$,}\\
               -nk(k+1)     &\text{otherwise,} \endcases
$$
$$
    \Null(L) = \Null_{-1}(L) =
        \cases 2k &\text{ if $n\equiv 0 \mod 4$,}\\
               0    &\text{otherwise.} \endcases
$$
For a proof I gave just a reference to [\refNeuI]. Now Proposition \propHirz\ provides
missing details of the computation because the link in question is $T_{0,1}(nk,2k)$.
\smallskip

\noindent{\bf Remark \remMultilink.}
If $m_1$ and $m_2$ are any non-negative integers, we define  $T_{m_1,m_2}(p,q)$
by the same formula as above but we interpret it as a multilink in the sense of [\refEN, \refNeuI].
In this case the equality $\sigma_\lambda^-=M^+_\theta-M^-_\theta$ and its proof hold true.

\demo{ Proof of Proposition \propHirz }

\midinsert
\centerline{\epsfxsize=30mm\epsfbox{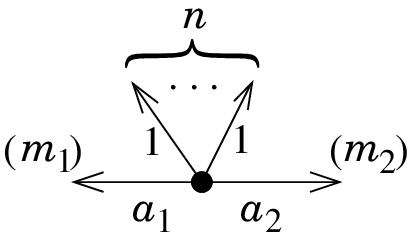}}
\botcaption{Figure \figSplicePQ} Splice diagram of $T_{m_1,m_2}(na_1,na_2)$
\endcaption
\endinsert

Let $p=a_1n$ and $q=a_2n$ where $n=\gcd(p,q)$. Then the splice diagram
of the (multi)link $T_{m_1,m_2}(p,q)$ is as shown in Figure~\figSplicePQ.
We have $M_0^+=n-1$, $M^-_0=0$, and the number of link components is
$n+m_1+m_2=(M_0^+ - M_0^-) + 1 + m_1+m_2$. It is also easy to check that
$$
    N_\varphi^+ - N_\varphi^- = (M_\varphi^+ - M_\varphi^-)
    + 2\sum_{0\le\theta<\varphi} (M_\theta^+ - M_\theta^-).
$$
Thus, by Corollary~\corNeumann, the computation of $\sigma_\lambda$
reduces to that of $\sigma_\lambda^-$.
To compute $\sigma_\lambda^-$, we apply Theorem~\thNeumann\
with $a_j=m_j=1$ for $j=3,\dots,n+2$. Let the $a'_j$, $b_j$, and $s_j$ be
defined as in \S\sectNeumann\ with the range of the indices appropriately changed.
We have $a'_1=a_2$, $a'_2=a_1$, and $a'_j=a_1a_2$ for $j\ge3$ whence
$$
     m = m_1 a_2 + m_2 a_1 + n a_1 a_2.
$$
Hence  the linear function $u$ involved in the definition of
$M_\theta^\pm$, $N_\theta^\pm$ takes the form
$u(x,y) = (a_2 x + a_1 y)/m$.
So, when $m\theta\not\in\Z$, we have $M_\theta^+ - M_\theta^- = 0 - 0 = \sigma_\theta^-$.
Let us fix $\theta = k/m$ with $k\in\Z$, $0<k<m$.
Then, by Theorem~\thNeumann, we have $\sigma_\theta^- = -2\sum\ld x_j\rd$
where $x_j = s_jk/m$, $j=1,\dots,n+2$.

We choose $b_1$ and $b_2$ so that $a_1b_2+a_2b_1=1$ and
we set $b_j=0$ for $j\ge 3$.
Then $\sum b_j a'_j=1$, hence
$$
   a_1\dots a_{n+2}\sum s_j = \sum(m_j - b_j m)a'_j = m - m\sum b_j a'_j = 0
$$
and we obtain $\sum x_j = \sum s_j k/m = 0$. We have also $0<x_j=k/m<1$ for $j\ge3$.
Thus we may apply Lemma \lemDede\ which yields
$$
    \sigma_\lambda^- = -n + 2\Card\big\{ l\in\Z\mid x_1 < l < -x_2\}
                        + \one_{\Z}(x_1) + \one_{\Z}(x_2).               \eqno(\eqPropHirz)
$$
The inequalities $x_1 < l < -x_2$ in (\eqPropHirz) can be transformed as follows:
$$
\xalignat5
      &x_1<l &&\Leftrightarrow & (m_1-b_1m)k/m &< a_1 l &&\Leftrightarrow &b_1 k + a_1 l &> m_1k/m,
\\
     &l<-x_2 &&\Leftrightarrow & (m_1-b_2m)k/m &<-a_2 l &&\Leftrightarrow &b_2 k - a_2 l &> m_2k/m.
\endxalignat
$$
Since $\varphi:t\mapsto(b_1 k + a_1 t, b_2 k - a_2 t)$ is a parametrization of the line $\{u=\theta\}$
such that $\varphi(\Z)=\Z^2\cap\{u=\theta\}$
and since $\varphi(m_1k/m)=A$ and $\varphi(m_2k/m)=B$ are the intersection points of $\{u=\theta\}$
with $\partial\Pi$ (see Figure \figParallelogramm), we can rewrite
(\eqPropHirz) in the form
$$
     \sigma_\lambda^-
      = -n + 2M_\theta^+ + \partial M_\theta^+
$$
where $\partial M_\theta^+$ stands for $\one_{\Z^2}(A) + \one_{\Z^2}(B)$.
By combining this equation with
$$
    n = M_0^- + \tfrac12\partial M_0^-
    = M_\theta^- + M_\theta^+ + \tfrac12(\partial M_\theta^- + \partial M_\theta^+)
    \quad\text{and}\quad \partial M_\theta^+ = \partial M_\theta^-,
                                                                     \eqno(\eqMM)
$$
we obtain the desired expression for $\sigma_\lambda^-$.

Now we apply Corollary~\corNeumann\ to compute $\Null_\lambda$. By (\eqAlex)
we have 
$$
   \Delta_L(t) = \frac{(t-1)(t^m-1)^{n+m_1+m_2-2}}
                  {(t^{m/a_1}-1)^{1-m_1}(t^{m/a_2}-1)^{1-m_2}}
$$
and by (\eqMM) we have $M_\theta^+ + M_\theta^- = n - \partial M_\theta^+$.
Then in the case $(m_1,m_2)=(1,1)$ we immediately conclude that $M^+_\theta + M^-_\theta$
is the multiplicity of $e^{2\pi i\theta}$ as a root of $\Delta_L(t)$.
We leave the cases $(m_1,m_2)=(1,0)$ and $(m_1,m_2)=(0,0)$ to the reader. \qed
\enddemo


\head\sectSeif. Multivariate signatures of Seifert links \endhead

In this section we compute the multivariate signatures and nullities of any Seifert
link with any orientations of its components. Proposition~\propColor\ combined
with Corollary~\corNeumann\ allows to compute $\sigma_L$ and $\Null_L$ for
any positive Seifert link $L=\seif(a_1,\dots,a_k;a_{k+1},\dots,a_n)$, $a_i\ge 0$.
Theorem~\thEN(2) combined with Proposition~\propSigRev\ allows to reduce the
general case to the case of positive Seifert links.
Note that a reversing of orientation of $\Sigma(a_1,\dots,a_n)$
(switching the color of the central node of a splice diagram)
changes the sign of $\sigma_L$.

In fact, in most cases the computation of $\sigma_L$ and $\Null_L$
for a positive Seifert link
immediately follows from the facts formulated in previous sections.
Indeed, let us consider the coordinates
$(\theta_1,\dots,\theta_n)$ on $(\sph^1\setminus\{1\})^n$, $u_j=e^{2\pi i\theta_j}$,
Then the restriction of $\sigma_L$ and $\Null_L$ to the diagonal
$\theta_1=\dots=\theta_n$ is determined by Theorem~\thNeumann\
and Proposition~\propSigSig. By (\eqAlex), the zero set of the multivariate
Alexander polynomial is a union of parallel hyperplanes transverse to the diagonal.
Thus the extension of $\sigma_L$ and $\Null_L$ from the diagonal
to $(\sph^1\setminus\{1\})^n$ is determined by Corollary~\corSR\ everywhere except
the hyperplanes corresponding to the multiple roots of $\Delta_L$.
Thus the only thing remaining to do is to extend $\sigma_L$ and $\Null_L$
to the multiple components of $\{\Delta_L=0\}$. This is done using the splicing formula
in Theorem~\thDFL\ combined with the observation that a Seifert link can be
spliced in many different ways. Let us give the exact statement and a formal proof.

\proclaim{ Proposition \propColor }
Let $L=L_1\cup\dots\cup L_k=\seif(a_1,\dots,a_k;\, a_{k+1},\dots,a_n)$
with non-negative $a_1,\dots,a_n$.
Let $u=(u_1,\dots,u_k)$ where
$u_j=\exp(2\pi\theta_j)$, $0<\theta_j<1$, 
Then $\sigma_L(u)$ and $\Null_L(u)$ depend only on the sum
$\sum_{j=1}^k a'_j\theta_j$ (the $a'_j$ are as in (\eqDefAprime)), thus
$$
    \sigma_L(u) = \sigma_L(\lambda,\dots,\lambda)
    = \sigma_\lambda(L) + \sum_{1\le i<j\le k}\lk(L_i,L_j)
$$
and
   $\Null_L(u) = \Null_L(\lambda,\dots,\lambda)=\Null_\lambda(L)$ for
$$
   \lambda=\exp(2\pi i\theta), \qquad
   \theta = \frac{a'_1\theta_1+\dots+a'_k\theta_k}{a'_1+\dots+a'_k}.
$$
\endproclaim

\demo{ Proof }
If $a_j=0$ for some $j$, then $(a_1,\dots,a_n)$ is a permutation of $(0,1,\dots,1)$,
hence either $L$ is a trivial $k$-component link
(if $j>k$; Theorem~\thEN(4)), or $L$ is the link in Figure \figHopf.
and then the result follows from [\refDFL;~Theorem~2.10].
So, we assume that all $a_j$ are positive.

If $n=2$, the statement is evident because $L$ is either an unknot or a
2-component Hopf link.
If $n=3$, the statement follows
from Corollary~\corSimpleRoot\ because
the multivariate Alexander polynomial is of the form
(see (\eqAlex))
$$
   \Delta_L(t_1,\dots,t_k) = f(t_1^{a'_1}\dots t_k^{a'_k}-1)
$$
where $f(t)$ is a polynomial without multiple roots, hence
the zero divisor of $\Delta_L$ considered as an analytic function of
$\theta_1,\dots,\theta_k$
is a union of parallel hyperplanes transverse to the main diagonal and
taken with multiplicity 1.

If $k=2$, we consider the splice in Figure~\figSpliceAA.
The obtained link coincides with $L$
by the edge contraction property (Theorem~\thEN(6)). Hence the result
follows from Theorem~\thDFL(b) combined with the statement of the lemma for $n=3$..

\midinsert
  \centerline{\epsfxsize=65mm\epsfbox{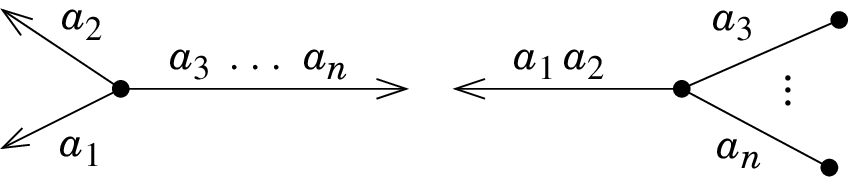}}
\botcaption{Figure \figSpliceAA} Splicing when $k=2$ in the proof of Prop.~\propColor.
\endcaption
\endinsert

So, we assume that $k\ge 3$.
We need to
prove that the signature and nullity
(considered as functions of $\theta_1,\dots,\theta_k)$ are
locally constant on each open $(k-1)$-dimensional polytope
$P_c=\{a'_1\theta_1+\dots+a'_k\theta_k=c\}\cap{]}0,1{[}^k$.
To this end it is enough to show that they are constant
on each open interval in $P_c$ defined by the condition that only two of the $\theta_k$'s
vary and the others are fixed. Moreover, it is enough to establish the constancy along
the subintervals of some fixed length.
Thus we reduce the problem to the following assertion.

\smallskip
{\sl Let $k\ge3$ and let $\tilde\theta_2$ and $\tilde\theta_3$ satisfy the conditions:
$$
   0 < \frac{\tilde\theta_2 - \theta_2}{a_2} =
    \frac{\theta_3 - \tilde\theta_3}{a_3} < \frac1{a_1}, \qquad
    0<\tilde\theta_2<1, \quad 0<\tilde\theta_3<1.                           \eqno(\eqCondA)
$$
Then
$\sigma_L(u) = \sigma_L(u_1,\tilde u_2,\tilde u_3,u_4,\dots,u_k)$
where $\tilde u_j=\exp(2\pi i\tilde\theta_j)$, $j=2,3$. }

\midinsert
\vbox{%
  \centerline{\lower-2ex\hbox{(i)  }\ \ \ \ \epsfxsize=68mm\epsfbox{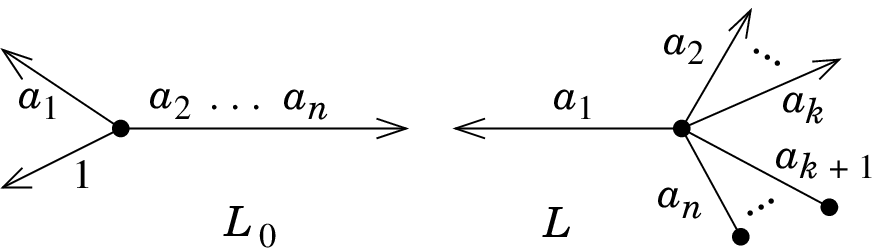}}
  \vskip 2ex
  \centerline{\lower-2ex\hbox{(ii) }\ \ \epsfxsize=70mm\epsfbox{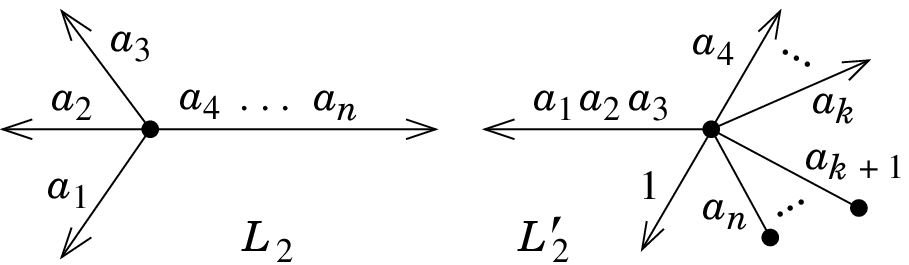}} }
\botcaption{Figure \figSpliceDE}
     Two splicings of $L_1$ in the proof of Prop.~\propColor.
\endcaption
\endinsert

\smallskip
Let us prove it.
Indeed, let $L_1$ be the link obtained by the splicing
according to Figure \figSpliceDE(i).
We have $L_1=\seif(1,a_1,\dots,a_k;\,a_{k+1},\dots,a_n)$
by the edge contraction property (Theorem~\thEN(6)).
By the same reason $L_1$ can be spliced also as in Figure~\figSpliceDE(ii).

We set $A_j=a_ja_{j+1}\dots a_n$
and $u_{\vec j}=(u_j,\dots,u_k)$,
in particular $u_{\vec1} = u$.
Due to the splicing in Figure \figSpliceDE(i), for any $u_0$
we have by Theorem~\thDFL(a):
$$
   \sigma_L(u) =
   \sigma_{L_1}(u_0,u^*_1,u_{\vec2})
   - \sigma_{L_0}(v,u_0,u^*_1)
   - \delta_L(u) \delta_{L_0}(v,u_0,u^*_1)                           \eqno(\eqSpliceD)
$$
where $u^*_1=u_1 u_0^{-a_1}$ and
$v = \prod_{j=2}^k u_j^{A_2/a_j}$
(the notation $\delta_L$ and $\delta_{L_0}$ must be clear).
This formula applies when $(u_1,v)\ne(1,1)$ which is the case because $u_1\ne 1$.

Similarly, if $(v_2,v_2')\ne (1,1)$, then the splicing in Figure~\figSpliceDE(ii) yields
$$
   \sigma_{L_1}(u_0,u^*_1,u_{\vec2}) =
   \sigma_{L_2}(v'_2,u^*_1,u_2,u_3) + \sigma_{L'_2}(v_2,u_0,u_{\vec4})
   + \delta_{L_2}(u^*_1,u_2,u_3) \delta_{L'_2}(u_0,u_{\vec4})    \eqno(\eqSpliceE)
$$
where $v'_2 = u_0^{A_4}\prod_{j=4}^k u_j^{A_4/a_j}$ and
$v_2 = (u^*_1)^{a_2 a_3} u_2^{a_1 a_3} u_3^{a_1 a_2}$.

Further, for the link $L_2$ (the left hand side of Figure~\figSpliceDE(ii)),
we consider the two splicings shown in Figure~\figSpliceFG. Then, if $(v_3,v'_3)\ne(1,1)\ne(v_4,v'_4)$,  we have

\vbox{
$$
   \sigma_{L_2}(v'_2,u^*_1, u_2, u_3)
    = \sigma_{L_3}(v'_3,u^*_1, u_2) + \sigma_{L'_3}(v_3, u_3,v'_2)
    + \delta_{L_3}(u^*_1, u_2)  \delta_{L'_3}(u_3,v'_2),           \eqno(\eqSpliceF)
$$
\vskip-30pt
$$
   \sigma_{L_2}(v'_2,u^*_1,\tilde u_2,\tilde u_3)
    = \sigma_{L_4}(v'_4,u^*_1,\tilde u_3) + \sigma_{L'_4}(v_4,\tilde u_2, v'_2)
    + \delta_{L_4}(u^*_1,\tilde u_3)  \delta_{L'_4}(\tilde u_2,v'_2)           \eqno(\eqSpliceG)
$$}
\noindent
where
$v_3  = (u^*_1)^{a_2} u_2^{a_1}$,
$v'_3 = (v_2')^{a_3} u_3^{A_4}$,
$v_4 = (u^*_1)^{a_3} \tilde u_3^{a_1}$,
$v'_4 = (v_2')^{a_2} \tilde u_2^{A_4}$.

\midinsert
\vbox{%
  \centerline{\lower-2ex\hbox{(i)  }\ \ \epsfxsize=70mm\epsfbox{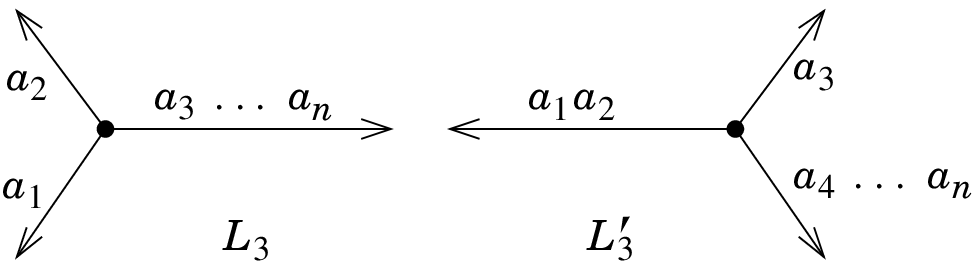}}\vskip 2ex
  \centerline{\lower-2ex\hbox{(ii) }\ \ \epsfxsize=70mm\epsfbox{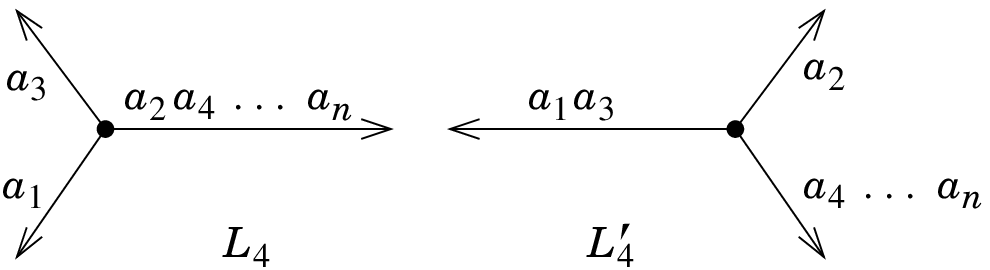}} }
\botcaption{Figure \figSpliceFG} Two splicings of $L_2$ in the proof of Prop.~\propColor.
\endcaption
\endinsert

Let $\theta^*\in\,{]}0,1{[}$ be such that $u^*_1=e^{2\pi i\theta^*_1}$
(recall that $u^*_1=u_1 u_0^{-a_1}$).
Let us choose $u_0$ so that $v_j\ne 1$ for $j=2,3,4$ (thus (\eqSpliceD)--(\eqSpliceG) hold) and $\theta^*_1$ satisfies the condition
$$
        0 < \tilde\theta_1 < \theta^*_1 < 1
       \quad\text{where $\tilde\theta_1$ is defined by}\quad
       \frac{\theta^*_1-\tilde\theta_1}{a_1} = \frac{\tilde\theta_2 - \theta_2}{a_2}
                                               = \frac{\theta_3 - \tilde\theta_3}{a_3}.   \eqno(\eqCondB)
$$
Such $u_0$ exists due to (\eqCondA).
Similarly to (\eqSpliceF) and (\eqSpliceG),
for $\tilde u_1=e^{2\pi i\tilde\theta_1}$ we have

\vbox{
$$
   \sigma_{L_2}(v'_2,\tilde u_1, \tilde u_2, u_3)
    = \sigma_{L_3}(v'_3,\tilde u_1,\tilde u_2) + \sigma_{L'_3}(v_3, u_3,v'_2)
    + \delta_{L_3}(\tilde u_1,\tilde u_2)  \delta_{L'_3}(u_3,v'_2),           \eqno(\eqSpliceFF)
$$
\vskip-30pt
$$
   \sigma_{L_2}(v'_2,\tilde u_1,\tilde u_2, u_3)
    = \sigma_{L_4}(v'_4,\tilde u_1, u_3) + \sigma_{L'_4}(v_4,\tilde u_2, v'_2)
    + \delta_{L_4}(\tilde u_1, u_3)  \delta_{L'_4}(\tilde u_2,v'_2).           \eqno(\eqSpliceGG)
$$}

\noindent
Since the statement of the lemma holds for $n=3$, the condition (\eqCondB) implies
$$
\split
    \sigma_{L_3}(v'_3,u^*_1, u_2) &= \sigma_{L_3}(v'_3,\tilde u_1,\tilde u_2),
    \qquad
    \delta_{L_3}(u^*_1, u_2) = \delta_{L_3}(\tilde u_1,\tilde u_2),
\\
    \sigma_{L_4}(v'_4,u^*_1,\tilde u_2) &= \sigma_{L_4}(v'_4,\tilde u_1, u_2),
    \qquad
    \delta_{L_4}(u^*_1,\tilde  u_2) = \delta_{L_4}(\tilde u_1,u_2).
\endsplit
$$
Hence all terms in the right hand sides of (\eqSpliceF) and (\eqSpliceG) are equal to the
corresponding terms in (\eqSpliceFF) and (\eqSpliceGG). Therefore the left hand sides are equal as well
and we obtain
$$
    \sigma_{L_2}(v'_2,        u_1^*,        u_2,        u_3) \overset{(\eqSpliceF),(\eqSpliceFF)}\to=
    \sigma_{L_2}(v'_2, \tilde u_1,   \tilde u_2,        u_3) \overset{(\eqSpliceG),(\eqSpliceGG)}\to=
    \sigma_{L_2}(v'_2,        u_1^*, \tilde u_2, \tilde u_3).
$$
Thus all terms in the right hand side of (\eqSpliceE) do not change when we replace
$u_2$ and $u_3$ with $\tilde u_2$ and $\tilde u_3$,
whence the same is true for (\eqSpliceD). This completes the proof for the signature.
The proof for the nullity is the same.
\qed\enddemo


\head\sectExample. Some examples
\endhead

\subhead\sectExample.1. Singularity link of two transverse cusps
\endsubhead
We start with the example considered in [\refNeuI, \S7].
Let $L$ be the link which is cut out by the curve
$$
      (x^3-y^2)(x^2-y^3)=0                 \eqno(\eqExample)
$$
on a sufficiently small sphere in $\C^2$ centered at the origin.
Then $L$ is obtained from the Hopf link by the $(2,3)$-cabling along each its component.
This corresponds to the splicing in Figure~\figExSplice\ where we write
the variable names near the corresponding arrowhead vertices.
We denote the splice components by $L_j\cup K_j$, $j=1,2$, where $L_j$ is
a trefoil knot and $K_j$ is its core.

\midinsert
\centerline{\epsfxsize=110mm\epsfbox{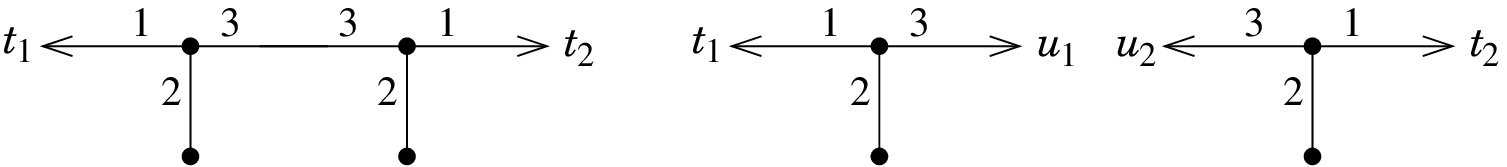}}
\botcaption{Figure \figExSplice}
    Splice diagram and splice components for
    the link in \S\sectExample.1
\endcaption
\endinsert

The Alexander polynomial of $L$ is
$$
   \Delta_L(t_1,t_2)=
   \frac{ (t_1^6 t_2^4 - 1)(t_1^4 t_2^6 - 1) }
        { (t_1^3 t_2^2 - 1)(t_1^2 t_2^3 - 1) }
   = (t_1^3 t_2^2 + 1)(t_1^2 t_2^3 + 1).
$$
and those of the splice components are $t_1^3u_1+1$ and $t_2^3u_2+1$.

\midinsert
\centerline{\epsfxsize=120mm\epsfbox{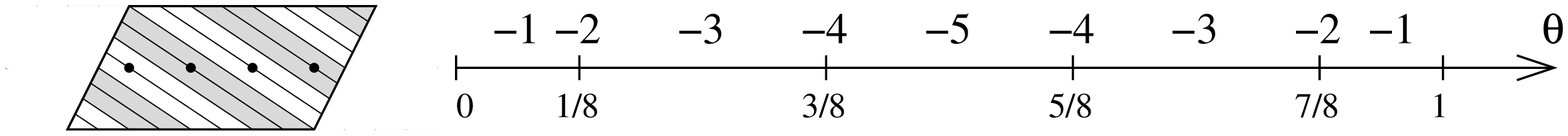}}
\botcaption{Figure \figExParal}
    $\Pi(3+1,0,1,2)$ and the values
    of $\sigma_\lambda(K_j\cup L_j)$, $\lambda=e^{2\pi i\theta}$
\endcaption
\endinsert

By Proposition~\propHirz, the Tristram--Levine signatures of each
splice component are
as in Figure~\figExParal.
Therefore, by Proposition~\propColor,
the multivariate signatures $\sigma_{K_j\cup L_j}(t_j,u_j)$ for $t_j=e^{2\pi i\theta_j}$,
$u_j=e^{2\pi i\varphi_j}$ are as in Figure~\figExSign\ on the left.
The nullity $\Null_{K_j\cup L_j}$ is $1-|\sign(\Delta_{L_j})|$ with only two exceptions:
$\Null_{K_j\cup L_j}(-1,1)=\Null_{K_j\cup L_j}(1,-1)=0$
(the small white circles in Figure~\figExSign).

By Theorem~\thDFL, for $(t_1,t_2)\in(\sph^1)^2$ we have
$$
   \sigma_L(t_1,t_2) = \sigma_{K_1\cup L_1}(t_2^2,t_1)
                     + \sigma_{K_2\cup L_2}(t_1^2,t_2)
                     + \delta_{(2)}(t_1)\delta_{(2)}(t_2)      \eqno(\eqDFL)
$$
unless $(t_1^2,t_2^2)=(1,1)$. In Figure~\figExSignRes\ we
show the values for two of the three terms of the right hand side of (\eqDFL)
(the pictures for $\sigma_{K_1\cup L_1}(t_2^2,t_1)$
and for $\sigma_{K_2\cup L_2}(t_1^2,t_2)$
are symmetric with respect to the diagonal $\theta_1=\theta_2$).
By summating the three terms we obtain the values of $\sigma_L$ shown in
Figure~\figExSignRes. We do not write the values of $\sigma_L$ on the line segments
composing the set $\Delta_L=0$, since by Lemma~\lemSR\ its value on each such segment
is the half-sum of the values on the two sides of the segments.

\midinsert
\centerline{ $\varphi_1$ \hskip 34mm $\theta_2$ \hskip34mm $\theta_2$ \hskip 22mm}
\vskip-8pt
\centerline{\epsfxsize=110mm\epsfbox{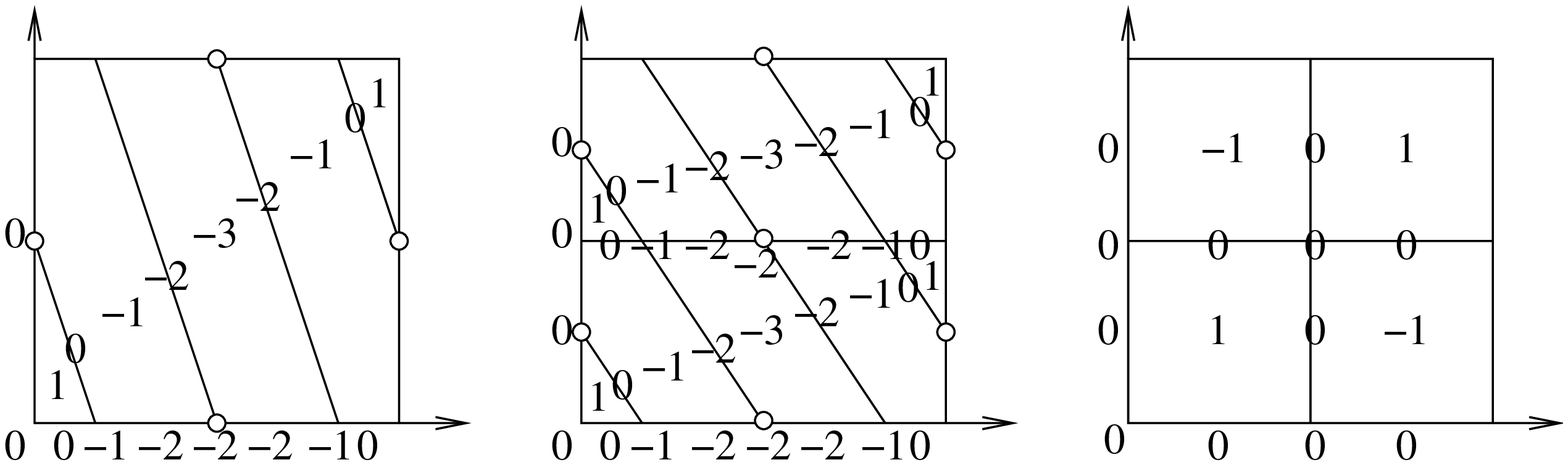}}
\vskip-25pt
\centerline{\hskip 30mm $\theta_1$ \hskip 34mm $\theta_1$ \hskip34mm $\theta_1$}
\vskip 17pt
\centerline{
    $\sigma_{K_1\cup L_1}(u_1,t_1)$ \hskip 15mm
    $\sigma_{K_1\cup L_1}(t_2^2,t_1)$ \hskip 15mm
    $\delta_{(2)}(t_1)\delta_{(2)}(t_2)$   }
\botcaption{Figure \figExSign}
    The values of summands in the right hand side of (\eqDFL)
\endcaption
\endinsert

In fact Theorem~\thDFL\ states that $\sigma_L$ is as shown in Figure~\figExSignRes\
only for $(t_1,t_2)\ne(\pm1,\pm1)$. However a correction term $\Delta\sigma$
can be computed using [\refDFLii, Thm.~5.3] and it appears to be zero.
In the following computation of $\Delta\sigma$ we use the terminology and
notation from [\refDFLii]. We perform the computation at $(t_1,t_2)=(-1,-1)$ only.
The Conway potential function of $K_j\cup L_j$ and that of the
trefoil knot $K_j$ are
$$
   \nabla_{K_j\cup L_j}(u_j,t_j) = t_j^3u_j + t_j^{-3}u_j^{-1},
\qquad
   \nabla_{K_j}(t_j) = (t_j^2-1+t_j^{-2})/(t-t^{-1})
$$
(see [\refCim]).
Hence the slopes are
$\kappa_j=(L_j/K_j)(-1)=-(\nabla_{K_j\cup L_j})'_{u_j}(1,i)/2\nabla_{K_j}(i)=2/3$,
$j=1,2$ (see [\refDFLii, Thm.~3.21]). Hence $\Delta\sigma=0$
(see [\refDFLii, Thm.~5.3 and Rem.~5.4]).

\midinsert
\centerline{\epsfxsize=60mm\epsfbox{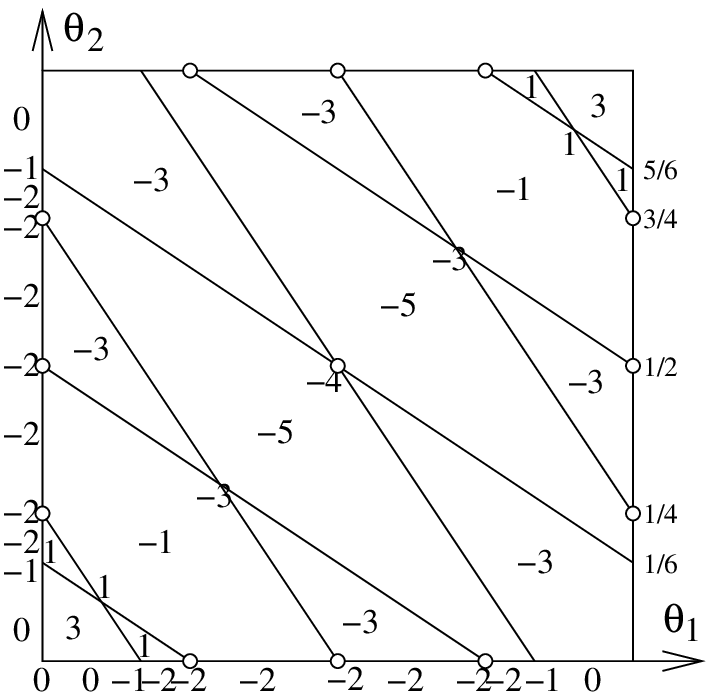}}
\botcaption{Figure \figExSignRes}
    The values of $\sigma_L(t_1,t_2)$, $t_j=e^{2\pi i\theta_j}$
\endcaption
\endinsert

Similarly to the signatures,
Theorem~\thDFL\ and [\refDFL, Thm.~5.3] also allow to compute
$\Null_L(t_1,t_2)$ for any $(t_1,t_2)\in(\sph^1)^2$. It is equal to the number of
components of $\Delta_L=0$ passing through the point $(t_1,t_2)$ unless
$(t_1,t_2)$ is one of $(-1,-1)$, $(1,i^k)$, or $(i^k,1)$, $k=1,2,3$
(the small white circles in Figure~\figExSignRes) where it is less by $1$.

By Proposition~\propSigSig, we deduce that the Tristram-Levine signatures of $L$
are as in Figure~\figExSignDiag\ and the non-zero nullities are:
$\Null_{-1}(L)=1$ and $\Null_\lambda(L)=2$ for $\lambda=e^{2\pi i k/10}$ with
$k\in\{1,3,7,9\}$. This fact well agrees with the Seifert form of $L$
computed in [\refNeuI, \S7] (see however Remarks~\remOrient\ and \remNeumann).

\midinsert
\centerline{\epsfxsize=100mm\epsfbox{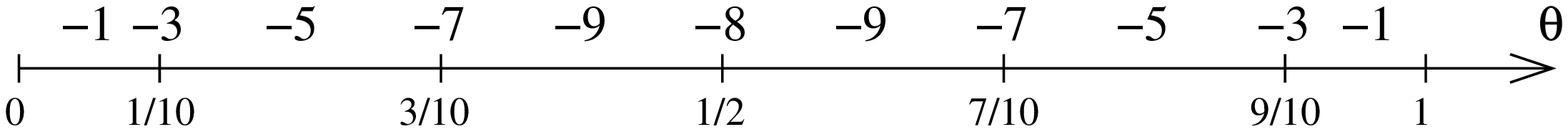}}
\botcaption{Figure \figExSignDiag}
    The values of $\sigma_\lambda(L)$, $\lambda=e^{2\pi i\theta}$
\endcaption
\endinsert

\medskip\noindent{\bf Question \sectExample.1.}
Is it true that the slopes of all iterated torus links
(all graph links in graph $3$-spheres?) are defined and can be computed
by [\refDFLii, Thm.~3.21]?
\medskip

\medskip
\subhead\sectExample.2. The link at infinity of two half-cubic parabolas
\endsubhead

Let $L^\infty$ be the link which is cut out by the same curve (\eqExample)
but on a sufficiently large sphere in $\C^2$ centered at the origin.
The splice diagram is obtained from Figure~\figExSplice\ by exchanging the labels
``$2$'' and ``$3$''. The Alexander polynomial is
$$
    (t_1^6t_2^4+t_1^3t_2^2+1)(t_1^4t_2^6+t_1^2t_2^3+1).
$$
Performing the same computations as in the previous example, we
obtain the values of $\sigma_{L^\infty}$ and $\Null_{L^\infty}$ as in
Figure~\figExTwo\ (as in \S\sectExample.1, the small white circles mark
the points where $\Null_{L^\infty}$ is one less than the number of
components of $\Delta_{L^\infty}$). The Tristram--Levine signatures of
$L^\infty$ are shown in Figure~\figExTwoDiag. The nullity is $1$ (resp. $2$)
at primitive 3rd (resp. $15$th) roots of unity.

\midinsert
\centerline{\epsfxsize=75mm\epsfbox{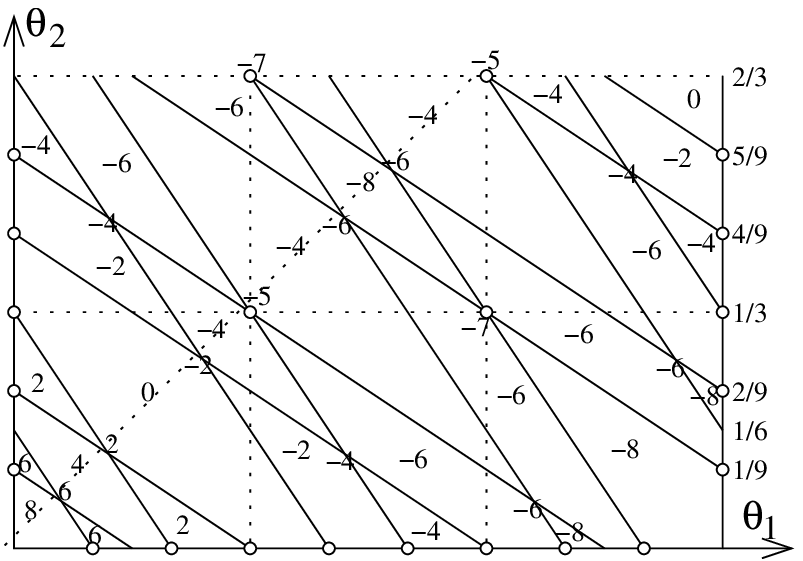}}
\botcaption{Figure \figExTwo}
    The values of $\sigma_{L^\infty}(t_1,t_2)$, $t_j=e^{2\pi i\theta_j}$
\endcaption
\endinsert

\midinsert
\centerline{\epsfxsize=125mm\epsfbox{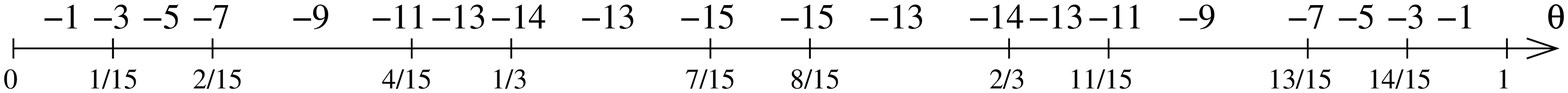}}
\botcaption{Figure \figExTwoDiag}
    The values of $\sigma_\lambda(L^\infty)$, $\lambda=e^{2\pi i\theta}$
\endcaption
\endinsert

\medskip\noindent{\bf Question \sectExample.2.}
One observes that the decomposition of the Hermitian isometric structure
of $L^\infty$ is given by the same formulas as in [\refNeuE], [\refNeuI]
with the only exception that the factors $\Lambda_\lambda^2$ appear with
the negative sign.
Is this true for all links at infinity of affine algebraic curves in $\C^2$
or for some natural subclass of such links?

\enddocument